\documentclass[draft,leqno, 10pt]{article}

\usepackage{amsmath,amsthm,amsfonts,amscd} 
\usepackage{eucal,eufrak}  
\usepackage{verbatim}     
\usepackage[dvips]{graphicx}

\newtheorem{thm}{Theorem}[section]
\newtheorem{cor}[thm]{Corollary}
\newtheorem{lemma}[thm]{Lemma}

\newcommand{\BB}[1]{\ensuremath{\mathbb{#1}}}
\newcommand{\Q}{\ensuremath{\BB{Q}}}
\newcommand{\R}{\ensuremath{\BB{R}}}

\newcommand{\C}{\ensuremath{\BB{C}}}

\DeclareMathOperator{\sgn}{sgn}
\DeclareMathOperator{\T}{T}
\DeclareMathOperator{\Res}{Res}
\DeclareMathOperator{\rec}{rec}

\numberwithin{equation}{section}

\begin{document}

\title{The distribution of Mahler's measures of reciprocal polynomials}
\author{Christopher D. Sinclair \\ The University of Texas at Austin \\
1 University Station C1200 \\ Austin, Texas 78757}
\maketitle

\begin{abstract}
We study the distribution of Mahler's measures of reciprocal
polynomials with complex coefficients and bounded even degree.  We
discover that the distribution function associated to Mahler's measure
restricted to monic reciprocal polynomials is a reciprocal (or
anti-reciprocal) Laurent polynomial on $[1,\infty)$ and identically
zero on $[0,1)$.  Moreover, the coefficients of this Laurent polynomial
are rational numbers times a power of $\pi$.  We are led to this
discovery by the computation of the Mellin transform of the
distribution function.  This Mellin transform is an even (or odd)
rational function with poles at small integers and residues that are
rational numbers times a power of $\pi$.  We also use this Mellin
transform to show that the volume of the set of reciprocal polynomials
with complex coefficients, bounded degree and Mahler's measure less
than or equal to one is a rational number times a power of $\pi$.
\end{abstract}

\vspace{.5cm}

\begin{section}{Introduction}
The Mahler's measure of a polynomial $f(x) \in \C[x]$ is given by the
expression 
\begin{equation}
\label{mahlers measure}
\mu(f) = \exp \left\{\int_0^1 \log | f(e^{2 \pi i t}) | \, dt \right\}.
\end{equation}
If $f(x)$ has degree $M$ and factors over $\C$ as $f(x) = w_M
\prod_{m=1}^M (x - \beta_n)$, then by Jensen's formula,
$$
\mu(f) = |w_M| \prod_{m=1}^M \max\{1, |\beta_m| \}.
$$ 
It is readily apparent that Mahler's measure is a
multiplicative function on $\C[x]$.  In this sense Mahler's measure
forms a natural height function on $\C[x]$.  In this paper we study
the distribution of values of Mahler's measure restricted to the set
of reciprocal polynomials with bounded even degree and complex
coefficients.

$f(x)$ is said to be reciprocal if it satisfies the condition
$$
x^M f\left(\frac{1}{x}\right) = f(x).
$$ 
If $f(x)$ is reciprocal and $f(x) = \sum_{m=0}^M w_m x^m$, then it
is easily seen that $w_m = w_{M-m}$ for $m = 0, \ldots, M$.  The
reciprocal condition also imposes a condition on the roots of $f(x)$:
if $f(\alpha) = 0$, then $f(\alpha^{-1}) = 0$.  If $M = 2N$ there
exists a Laurent polynomial
\begin{equation}
\label{reciprocal polynomial}
p_{\mathbf{v}}(x) = v_0 + \sum_{n=1}^N v_n \left(x^n +
x^{-n}\right)
\end{equation}
such that $f(x) = x^N p_{\mathbf{v}}(x)$.  We call $p_{\mathbf{v}}(x)$
the reciprocal Laurent polynomial with coefficient vector 
$\mathbf{v}$.  The collection of reciprocal Laurent polynomials with
complex coefficients forms a graded algebra.  

The integral defining Mahler's measure makes sense for reciprocal
Laurent polynomials, and it is easily seen that $\mu(p_{\mathbf{v}}) = 
\mu(f)$.  It is convenient to work with reciprocal Laurent
polynomials since they form an algebra (the set of reciprocal
polynomials is not closed under addition).  We define the {\it
reciprocal Mahler's measure} to be the function $\mu_{\rec}:\C^{N+1}
\rightarrow \R$ given by
\begin{equation}
\label{reciprocal mahlers measure}
\mu_{\rec}(\mathbf{v}) = \mu(p_{\mathbf{v}}) =
\exp\left\{\int_0^1 \log\left|v_0 + 
2 \sum_{n=1}^N v_n \cos(2 \pi n t) \right| \, dt \right\}.
\end{equation}

If $\mathbf{v} = (v_0, \ldots, v_L, 0, \ldots, 0)$ with $v_L \neq 0$,
then there exist $\alpha_1, \ldots, \alpha_{2L}$ not necessarily
distinct nonzero complex roots of $p_{\mathbf{v}}(x)$.  By reordering
if necessary, we may assume $\alpha_{L+n} = \alpha_n^{-1}$, and we may
write
\begin{equation*}
\label{root polynomial}
x^L p_{\mathbf{v}}(x) = v_L \prod_{n = 1}^L (x - \alpha_n)(x -
\alpha_n^{-1}),
\end{equation*}
and from Jensen's formula we have
\begin{equation}
\label{jensens formula}
\mu_{\rec}(\mathbf{v}) = |v_L| \prod_{n=1}^L \max\{|\alpha_n|, \left|
\alpha_n^{-1}\right|\}.
\end{equation}
From this expression we see for all $\mathbf{v} \in \C^{N+1}$ and $k
\in \C$ the reciprocal Mahler's measure is: 
\begin{enumerate}
\item \label{non-negative property} non-negative:
$\mu_{\rec}(\mathbf{v}) \geq 0$,
\item \label{homogeneous property} homogeneous: $\mu_{\rec}(k \,
  \mathbf{v}) = |k|\,\mu_{\rec}(\mathbf{v})$, and
\item \label{positive-definite property} positive-definite:
  $\mu_{\rec}(\mathbf{v}) = 0 \quad \mbox{if and only if} \quad
  \mathbf{v} = \boldsymbol{0}$.
\end{enumerate}
In addition $\mu_{\rec}$ is continuous as originally proved by Mahler
\cite{mahler}.

By properties (i), (ii) and continuity, we find that $\mu_{\rec}$
is a symmetric distance function in the sense
of the geometry of numbers (see for instance the discussion in
\cite[chapter IV]{cassels}).  $\mu_{\rec}$ satisfies all the
properties of a metric except the triangle inequality.  The `unit
ball' is thus not convex.  Explicitly,
$$
\mathcal{V}_{N+1} = \{\mathbf{v} \in \C^{N+1}\, : \,
\mu_{\rec}(\mathbf{v}) \leq 1\} 
$$ 
is a symmetric star body.  By property (iii) this star body is
bounded.  We call $\mathcal{V}_{N+1}$ the complex star body determined
by the reciprocal Mahler's measure.  One of the principal results
presented here is the computation of the volume (Lebesgue measure) of
$\mathcal{V}_{N+1}$.

We introduce the {\it monic} reciprocal Mahler's measure, $\nu_{\rec}:
\C^N \rightarrow \R$, defined by
$$
\nu_{\rec}(\mathbf{b}) = \mu_{\rec} \left( 
    \begin{array}{c}  \mathbf{b} \\ 1
    \end{array} 
\right).
$$
Thus $\nu_{\rec}(\mathbf{b})$ is the Mahler's measure of the monic reciprocal
Laurent polynomial
\begin{equation}
\label{monic reciprocal laurent polynomial}
\tilde{p}_\mathbf{b}(x) = \left(x^{N} + x^{-N}\right)+ b_0 +
\sum_{n=1}^{N-1} b_n \left(x^n + x^{-n}\right).
\end{equation}
We denote Lebesgue measure on Borel subsets of $\C^N$ by
$\lambda_{2N}$, and introduce the distribution function associated with
the monic reciprocal Mahler's measure, $h_N(\xi):[0, \infty)
\rightarrow [0,\infty)$, given by
\begin{equation*}
h_N(\xi) = \lambda_{2N}\left\{ \mathbf{b} \in \C^N :
\nu_{\rec}(\mathbf{b}) \leq \xi \right\}.
\end{equation*}
$h_N(\xi)$ encodes statistical information about the distribution of
Mahler's measures of reciprocal polynomials with complex coefficients
and even degree bounded by $2N$.  

The distribution function $h_N(\xi)$ is increasing and continuous from
the right.  By equation \ref{jensens formula} $\nu_{\rec}(\mathbf{b})
\geq 1$ for all $\mathbf{b} \in \C^N$, thus $h_N(\xi)$ is identically
zero on $[0, 1)$.  In fact, $h_N(1) = 0$.  To see this, suppose
$\mathbf{b} \in \C^N$ with $\nu_{\rec}(\mathbf{b}) = 1$.  Then, from
equation \ref{jensens formula}, $\tilde{p}_{\mathbf{b}}(x)$ has all
its roots on the unit circle.  Thus, if $\alpha$ is a root of
$\tilde{p}_{\mathbf{b}}(x)$ then so is $\overline{\alpha} =
\alpha^{-1}$.  We find that $\mathbf{b} \in \R^N$, and hence the set
of $\mathbf{b} \in \C^N$ such that $\nu_{\rec}(\mathbf{b}) = 1$ has
$\lambda_{2N}$-measure 0.  Thus $h_N(1) = 0$, and $h_N(\xi)$ is
continuous at $\xi = 1$.

We recall the definition of the Mellin transform.  Given a function
$g: [0,\infty) \rightarrow \R$, the Mellin transform of $g$ is the
function of the complex variable $s$ given by
$$
\widehat{g}(s) = \int_0^{\infty} \xi^{-2s} g(\xi) \, \frac{d \xi}{\xi}.
$$

We will give an explicit formula for $h_N(\xi)$ by computing its
Mellin transform.  We note that, since $h_N(\xi)$ is identically zero
on $[0,1]$ the integral defining $\widehat{h_N}(s)$ can be written
with domain of integration $[1, \infty)$.

The integral defining $\widehat{h_N}(s)$ converges in the half plane
$\Re(s) > N$.  To see this, we use the following consequence of
Jensen's inequality: 
$$
\mu(f) \leq \| f(x) \|_2,
$$
where $\|f(x) \|_2$ is the Euclidean norm of the coefficient vector
of $f(x)$.  Thus from equation \ref{monic reciprocal laurent
polynomial} we have
$$\nu_{\rec}(\mathbf{b}) \leq (2 + |b_0|^2 + 2|b_1|^2 + \ldots +
2|b_{N-1}|^2)^{1/2} \leq \sqrt{2}(1 + |b_0|^2 + \ldots +
|b_{N-1}|^2)^{1/2},$$ 
and hence
$$
\{\mathbf{b} \in \C^{N} : \nu_{\rec}(\mathbf{b}) \leq \xi\}
\subset
\left\{\mathbf{b} \in \C^{N} : (1 + |b_0|^2 + \ldots + |b_{N-1}|^2)^{1/2}
\leq \frac{\xi}{\sqrt{2}}\right\}.
$$ 
The latter set is a `slice' of a solid sphere of dimension $2N + 1$,
and is thus a solid sphere of dimension $2N$ with radius less than
$\xi/\sqrt{2}$.  Thus there exists a constant $C$ such that
$$
h_N(\xi) = \lambda_{2N}\{\mathbf{b} \in \C^{N} :
\nu_{\rec}(\mathbf{b}) \leq \xi\} \leq C \xi^{2N}.
$$
It follows that
$$
\widehat{h_N}(s) = \int_1^{\infty} \xi^{-2s} h_N(\xi) \, \frac{d\xi}{\xi} 
\leq C \int_1^{\infty} \xi^{2N - 2s} \, \frac{d\xi}{\xi}.
$$
The latter integral converges if $\Re(s) > N$, and hence
$\widehat{h_N}(s)$ is defined in the half plane $\Re(s) > N$.

We follow the method introduced by Chern and Vaaler in
\cite{chern-vaaler} to express the volume of $\mathcal{V}_{N+1}$ in
terms of the Mellin transform of $h_N(\xi)$. 
\begin{thm}
\label{volume theorem}
For each positive integer $N$
$$\lambda_{2N+2}(\mathcal{V}_{N+1}) = 2\pi \widehat{h_N}(N+1).$$
\end{thm}
\begin{proof}
The volume of
$\mathcal{V}_{N+1}$ is given by
\begin{equation}
\label{volume integral}
\lambda_{2N+2}(\mathcal{V}_{N+1}) = \int_{\C} \lambda_{2N}\left\{
\mathbf{b} \in \C^N : \mu_{\rec} \left( 
    \begin{array}{c}  \mathbf{b} \\ z
    \end{array} 
\right) \leq 1 \right\} \, d\lambda_2(z).
\end{equation}
By the homogeneity of $\mu_{\rec}$ we see
\begin{eqnarray*}
\lambda_{2N}\left\{
\mathbf{b} \in \C^N : \mu_{\rec} \left( 
    \begin{array}{c}  \mathbf{b} \\ z
    \end{array} 
\right) \leq 1 \right\}
&=&
\lambda_{2N}\left\{
z\mathbf{c} \in \C^N : \mu_{\rec} \left( 
    \begin{array}{c}  z \mathbf{c} \\ z
    \end{array} 
\right) \leq 1 \right\}\\
&=&
|z|^{2N}\lambda_{2N}\left\{
\mathbf{c} \in \C^N : \mu_{\rec} \left( 
    \begin{array}{c} \mathbf{c} \\ 1
    \end{array} 
\right) \leq \frac{1}{|z|} \right\} \\
&=& |z|^{2N} h_N \left( \frac{1}{|z|} \right)
\end{eqnarray*}
and thus the integral in equation \ref{volume integral} can be written as
$$
\int_{\C} |z|^{2N} h_N\left(\frac{1}{|z|}\right) \, d\lambda_2(z) \\
=2\pi \int_0^1 r^{2N+1} h_N\left(\frac{1}{r}\right) \, dr.
$$
The domain of integration in the latter integral is $[0,1)$ since
$h_N(1/r)$ is identically zero on $[1,\infty)$.  By the change of variables
$r=1/ \xi$ we find 
$$
\lambda_{2N+2}\left(\mathcal{V}_{N+1}\right)
=2\pi \int_1^\infty {\xi}^{-2(N+1)-1} h_N(\xi) \, d\xi
= 2\pi \widehat{h_N}(N+1).
$$
\end{proof}

If we regard the integral defining $\widehat{h_N}(s)$ as a
Lebesgue-Stieltjes integral, we may use integration by parts to 
write
\begin{equation*}
\label{integration by parts}
\widehat{h_N}(s) = \left. -\frac{\xi^{-2s} h_N(\xi)}{2s}
\right|_1^{\infty} + \frac{1}{2s} \int_1^{\infty} \xi^{-2s} 
\,dh_N(\xi).
\end{equation*}
Since $h_N(1) = 0$ and $h_N(\xi)$ is dominated by $C \xi^{2N}$, the
first term vanishes when $\Re(s) > N$. After a change of variables,
we can write
\begin{equation*}
\label{nu integral}
\widehat{h_N}(s) = \frac{1}{2s} \int_{\C^N} \nu_{\rec}(\mathbf{b})^{-2s} \,
d\lambda_{2N}(\mathbf{b}).
\end{equation*}
The latter integral is interesting enough to name:
$$H_N(s) = \int_{\C^N} \nu_{\rec}(\mathbf{b})^{-2s} \,
d\lambda_{2N}(\mathbf{b}).$$

The bulk of this paper is committed to the discovery that
$H_N(s)$ analytically continues to a rational function of $s$.
\begin{thm}
\label{main theorem}
For each positive integer $N$, the function $H_N(s)$ extends by
analytic continuation to an (even or odd) rational function.  In
particular,
\begin{equation*}
\label{rational function}
H_N(s) = \prod_{n=1}^N \frac{2 \pi s}{s^2 - n^2}.
\end{equation*}
\end{thm}

\begin{cor}
For each positive integer $N$,
$$\lambda_{2N+2}(\mathcal{V}_{N+1}) = \frac{2^N \, \pi^{N+1} \,
    (N+1)^N}{(2N + 1)!}.$$
\end{cor}
\begin{proof}
This follows immediately from Theorem \ref{volume theorem} and Theorem
\ref{main theorem}.
\end{proof}
\begin{cor}
\label{mellin inversion corollary}
For each positive integer $N$, $h_N(\xi)$ is a reciprocal or
anti-reciprocal Laurent polynomial on the domain $[1,\infty)$ and
identically zero on $[0, 1)$.  Explicity, if $\xi \geq 1$ then
$$h_N(\xi) = 2^N \pi^N \sum_{n = 1}^N \frac{(-1)^{N-n}\,
    n^N}{(N+n)!(N-n)!} \left(\xi^{-2 n} + (-1)^N \xi^{2 n} \right)
$$
\end{cor}
\begin{proof}
$\widehat{h_N}(s) = H_N(s)/2s$ is a rational function whose
denominator is a product of distinct linear factors of the form $s-n$.
We use partial fraction decomposition to write 
$$
\widehat{h_N}(s) = \sum_{n=1}^N \left( \frac{\rho(n)}{s - n} +
\frac{\rho(-n)}{s + n} \right) \quad \mbox{where} \quad
\rho(n) = \mathop{\Res}_{s = n}(\widehat{h_N}(s)).
$$
We compute $\rho(n)$:
$$(s - n)\widehat{h_N}(s) = \frac{\pi}{s + n} \prod_{\scriptstyle m = 1 \atop \scriptstyle m \neq n}^N \frac{2 \pi s}{s^2 - m^2}$$
and so
\begin{eqnarray}
\rho(n) &=& \pi^N 2^{N-2} n^{N-2} \prod_{m=1}^{n-1} \frac{1}{(n^2-m^2)}
\prod_{m = n+1}^N \frac{1}{(n^2-m^2)} \nonumber \\
&=& \pi^N 2^{N-2}n^{N-2}
\left(\frac{n!}{(n-1)!(2n-1)!}\right)\left(\frac{(-1)^{N-n}
  (2n)!}{(N+n)!(N-n)!}\right) \nonumber \\
\label{residue}
&=& \pi^N 2^{N-1} n^{N} (-1)^{N-n} \frac{1}{(N+n)!(N-n)!}.
\end{eqnarray}
It is clear that $\rho(-n) = (-1)^{N} \rho(n)$, and so
$$
\widehat{h_N}(s) = \sum_{n=1}^N \rho(n) \left( \frac{1}{s - n} +
\frac{(-1)^{N}}{s + n} \right)
$$
A quick calculation shows, for $s > n$,
$$
2 \int_1^{\infty} \xi^{-2s} (\xi^{-2n} \pm \xi^{2n}) \, \frac{d\xi}{\xi}
= \frac{1}{s + n} \pm \frac{1}{s - n}.
$$
And so, by the uniqueness of the Mellin transform we find 
\begin{equation}
\label{h_N Laurent polynomial}
h_N(\xi) = \sum_{n=1}^N 2 \rho(n) (\xi^{-2n} + (-1)^N \xi^{2n})
\end{equation}
for $\xi \in (1, \infty)$.  The lemma follows by substituting
equation \ref{residue} into equation \ref{h_N Laurent
polynomial}.
\end{proof}

We outline the proof of Theorem \ref{main theorem}.  Given
$\boldsymbol{\alpha} \in (\C \setminus \{0\})^N$, we can create the
unique monic reciprocal Laurent polynomial $\tilde{p}_{\mathbf{a}}(x)$ having
$\alpha_1, \ldots, \alpha_N, \alpha_1^{-1}, \ldots, \alpha_N^{-1}$ as
roots.  We will use the change of variables $\boldsymbol{\alpha}
\mapsto \mathbf{a}$ to write $H_N(s)$ as an integral over root vectors
of reciprocal Laurent polynomials, as opposed to coefficient vectors.
This change of variables is useful, since by equation \ref{jensens
formula}, $\nu_{\rec}(\mathbf{a})$ is a simple product in the roots of
$\tilde{p}_\mathbf{a}(x)$ ({\it i.e.} in the coordinates of
$\boldsymbol{\alpha}$).  Analysis of the Jacobian of this change of
variables will allow us to write $H_N(s)$ as the determinant of an $N
\times N$ matrix, the entries of which are Mellin transforms which
evaluate to rational functions of $s$.  Theorem \ref{main theorem}
will follow from the evaluation of the determinant of this matrix.

Before proceding to the proof of Theorem \ref{main theorem}, we present
$\mu_{\rec}$ and $\mathcal{V}_{N+1}$ from another perspective.  Given
the positive integer $M$, we define the Mahler's measure function to be 
$\mu: \C^{M+1} \rightarrow \R$ where $\mu(\mathbf{u})$ is the
Mahler's measure of the polynomial with coefficient vector
$\mathbf{u}$.  As was shown in $\cite{chern-vaaler}$, $\mu$ is
non-negative, homogeneous, positive-definite and continuous.  Thus
$\mu$ is a symmetric distance function and the set
$$
\mathcal{U}_{M+1} = \left\{\mathbf{u} \in \C^{M+1} : \mu(\mathbf{u})
\leq 1\right\}
$$ 
is a bounded symmetric starbody.  Let $M = 2N$ and consider the
linear map $\Lambda: \C^{N+1} \rightarrow \C^{2N+1}$ defined by
$\Lambda(\mathbf{v}) = \left(v_0, v_1, \ldots, v_{N-1}, v_N, v_{N-1},
\ldots, v_1, v_0\right)^{\T}$.  We define $V = \Lambda(\C^{N+1})$ to
be the subspace of reciprocal coefficient vectors.  By equations
\ref{mahlers measure}, \ref{reciprocal polynomial} and \ref{reciprocal
mahlers measure} we find $\mu_{\rec}(\mathbf{v}) =
\mu(\Lambda(\mathbf{v}))$. Thus, the starbody formed by the
intersection $\mathcal{U}_{2N+1}$ and $V$ is related to the reciprocal
starbody.  Specifically,
$$
\mathcal{V}_{N+1} = \Lambda^{-1}\left(V \cap \mathcal{U}_{2N+1}\right).
$$
Every bounded symmetric starbody uniquely determines a symmetric
distance function \cite[Chapter IV.2 Theorem 1]{cassels}.  Thus, armed
with $\mu$ and $\Lambda$, we could `discover' $\mu_{\rec}$.  Equation
\ref{jensens formula} can be recovered from the symmetry in the
definition of $\Lambda$, so we would lose no information if we were
to define $\mu_{\rec}$ in this manner.  

The volume of $\mathcal{U}_{M+1}$, as well as the subspace volume of the
starbody formed by intersecting $\mathcal{U}_{M+1}$ with the subspace
of real coefficient vectors was investigated in \cite{chern-vaaler}.
Thus the computation of the volume of $\mathcal{V}_{N+1}$ yields
subspace volume information of another `slice' of $\mathcal{U}_{2N+1}$.
\end{section}

\begin{section}{A change of variables}

Let $\C^{\times} = \C \setminus \{0\}$, and define the map
$\mathcal{E}_N: (\C^{\times})^N \rightarrow \C^N$ by
$\mathcal{E}_N(\boldsymbol{\alpha}) = \mathbf{a}$, where
$$
x^N \tilde{p}_\mathbf{a}(x) = \prod_{n = 1}^N (x + \alpha_n)(x +\alpha_n^{-1}).
$$ Thus the $n$th coordinate function of $\mathcal{E}_N$ is given by
$\varepsilon_n(\alpha_1, \ldots, \alpha_N, \alpha_1^{-1},
\ldots,\alpha_N^{-1})$, where $\varepsilon_n$ is the $n$th elementary
symmetric function in $2N$ variables.  Let $E_N: \C^N \rightarrow
\C^N$ be the function whose $n$th coordinate function is $e_n$, the
$n$th elementary symmetric function in $N$ variables.  That is, given
$\boldsymbol{\beta} \in \C^N$, if $\mathbf{b} =
E_N(\boldsymbol{\beta})$ then
$$
\prod_{n=1}^N (x + \beta_n) = x^N + \sum_{n=0}^{N-1} b_n x^n.
$$
It is well known that the Jacobian of $E_N(\boldsymbol{\beta})$ is
given by $|V(\boldsymbol{\beta})|^2$, where
\begin{equation}
\label{vandermonde}
V(\boldsymbol{\beta}) = \prod_{1 \leq m < n \leq N}(\beta_n - \beta_m)
=
 \det
\begin{pmatrix}
1         &    1      & \ldots & 1         \\
\beta_1   & \beta_2   & \ldots & \beta_N   \\
\beta_1^2 & \beta_2^2 & \ldots & \beta_N^2 \\
\vdots    & \vdots    & \ddots & \vdots    \\
\beta_1^{N-1} & \beta_2^{N-1} & \ldots & \beta_N^{N-1} 
\end{pmatrix}
\end{equation}
is the Vandermonde determinant.
We will relate the Jacobian of $\mathcal{E}_N$ to the Jacobian of
$E_N$. 

\begin{lemma}
\label{jacobian corollary}
For each positive integer $N$, the Jacobian of
$\mathcal{E}_N(\boldsymbol{\alpha})$ is given by 
$$
\left| V\left(\alpha_1 + \frac{1}{\alpha_1}, \ldots , \alpha_N +
\frac{1}{\alpha_N}\right) \right|^2 
\cdot
\prod_{n=1}^N \left| \left(\frac{\alpha_n^2 - 1}{\alpha_n^2}
\right)\right|^2. 
$$
\end{lemma}
\begin{proof}
By definition $\varepsilon_n(x_1, \ldots, x_N, x'_1, \ldots, x'_N)$ is
composed of all monomials of degree $n$ in the variables $x_1, \ldots
x_N, x'_1, \ldots x'_N$.  If we impose the relation $x_m x'_m = 1$ for
$m=1, \ldots, N$, then $\varepsilon_n(x_1, \ldots, x_N, x'_1, \ldots,
x'_N)$ is no longer homogeneous.  In this situation it is easy to see
that the monomials of degree $n$ of $\varepsilon_n(x_1, \ldots, x_N,
x'_1, \ldots, x'_N)$ consist of monomials which do not contain both
$x_m$ and $x'_m$ for $m=1, \ldots, N$.  Hence,
$$
\varepsilon_n(x_1, \ldots, x_N, x'_1, \ldots, x'_N) = e_n(x_1 + x'_1,
\ldots, x_N + x'_N) + (\mbox{monomials of degree } < n).
$$ 
In general $\varepsilon_n(x_1, \ldots, x_N, x'_1, \ldots, x'_N)$
contains monomials of degree $n-2M$ formed from monomials which
contain $x_m$ and $x'_m$ where $m$ runs over a subset of $1, \ldots,
N$ of cardinality $M$.  By counting the number of times each monomial
of degree $n-2M$ appears we arrive at the following identity. 
\begin{eqnarray*}
\varepsilon_n\left(\alpha_1, \ldots, \alpha_N, \frac{1}{\alpha_1}, \ldots,
\frac{1}{\alpha_N}\right) &=& e_n\left(\alpha_1 + \frac{1}{\alpha_1},
\ldots, \alpha_N + \frac{1}{\alpha_N}\right) \\ &+& {N - n - 2 \choose 1}
e_{n-2}\left(\alpha_1 + \frac{1}{\alpha_1}, \ldots, \alpha_N +
\frac{1}{\alpha_N}\right) \\ &+& {N -n - 4 \choose 2} e_{n-4}\left(\alpha_1 +
\frac{1}{\alpha_1}, \ldots, \alpha_N +\frac{1}{\alpha_N}\right) +
\ldots \\
&=& \sum_{M=0}^{[N/2]} {N -n - 2M \choose M} e_{n-2M}\left(\alpha_1 +
\frac{1}{\alpha_1}, \ldots, \alpha_N +\frac{1}{\alpha_N}\right) 
\end{eqnarray*}
where $[N/2]$ is the integer part of $N/2$.

Thus
$$
\mathcal{E}_N(\boldsymbol{\alpha}) = 
\begin{pmatrix}
1 & 0   & \ldots     & 0    \\
\ast      & 1 & \ldots     & 0    \\
\vdots & \vdots & \ddots   & \vdots\\
\ast      & \ast      & \ldots      & 1
\end{pmatrix}
E_N(\boldsymbol{\beta}).
$$ where 
$$\boldsymbol{\beta} = \left(\alpha_1 + \frac{1}{\alpha_1}, \ldots ,
\alpha_N + \frac{1}{\alpha_N}\right)^{\T}$$ and $\ast$ represents
entries which are not necessarily 0.  The Jacobian of
$E_N(\boldsymbol{\beta}) = |V(\boldsymbol{\beta})|^2$, and thus by
the chain rule we arrive at the formula for the Jacobian of
$\mathcal{E}_N(\boldsymbol{\alpha})$ given in the statement of the
lemma.
\end{proof}
The Jacobian of $\mathcal{E}_N(\boldsymbol{\alpha})$ is nonzero for
$\lambda_{2N}$-almost all points of $(\C^{\times})^N$, and
there are $2^N N!$ preimages for $\lambda_{2N}$-almost all
$\mathbf{a} \in \C^N$.  Employing the change of variables formula, we
find 
\begin{eqnarray}
H_N(s) &=& \int_{\C^N} \nu_{\rec}(\mathbf{a})^{-2s} \,
d\lambda_{2N}(\mathbf{a}) \nonumber \\
 &=& \frac{1}{2^N N!}\int_{(\C^{\times})^N}\left( \prod_{n=1}^N
\max\left\{|\alpha_n|, |\alpha_n^{-1}|\right\}^{-2s}
\left|\left(\frac{\alpha_n^2 - 1}{\alpha_n^2}
\right)\right|^2 \right) \nonumber  \\ 
\label{H_N(s) after jacobian}
&& \hspace{2cm} \times 
\left| V\left(\alpha_1 + \frac{1}{\alpha_1}, \ldots ,
\alpha_n + \frac{1}{\alpha_n}\right) \right|^2 
d\lambda_{2N}(\boldsymbol{\alpha}).
\end{eqnarray}
The latter integral admittedly looks formidable, however this change
of variables is beneficial since $\nu_{\rec}$ is a simple product.
\end{section}

\begin{section}{$H_N(s)$ is a determinant}
We first prove a short technical lemma concerning determinants.
\begin{lemma}
\label{determinant lemma}
Let $N$ be a positive integer. If $I = I(j, k)$ is an $N \times N$
matrix and $S_N$ is the $N$th symmetric group, then
\begin{equation}
\label{determinant lemma equation}
\det(I) =
\frac{1}{N!} \sum_{\tau \in S_N} \sum_{\sigma \in S_N} \sgn(\tau)
  \sgn(\sigma) \prod_{n=1}^N I(\tau(n), \sigma(n)).
\end{equation}
\end{lemma}
\begin{proof}
$$
\prod_{n=1}^N I\left(\tau(n), \sigma(n)\right) = \prod_{n=1}^N I(n, \sigma
\circ \tau^{-1}(n)).
$$ Thus we can write (\ref{determinant lemma
  equation}) as: 
\begin{eqnarray*}
\frac{1}{N!}\sum_{\tau \in S_N} \sum_{\sigma \in S_N} \sgn(\sigma
\circ \tau^{-1}) \prod_{n=1}^N I(n, \sigma \circ \tau^{-1}(n)) 
&=& \frac{1}{N!}\sum_{\tau \in S_N} \sum_{\sigma \in S_N} \sgn(\sigma)
\prod_{n=1}^N I(n, \sigma(n)) \\
&=& \sum_{\sigma \in S_N} \sgn(\sigma) \prod_{n=1}^N I(n, \sigma(n)),
\end{eqnarray*}
which is the familiar formula for $\det(I)$.
\end{proof}

Using equation \ref{vandermonde}, we expand the Vandermonde determinant
as a sum over the symmetric group to find
$$
\left| V\left(\alpha_1 + \frac{1}{\alpha_1}, \ldots ,
\alpha_n + \frac{1}{\alpha_n}\right) \right|^2 
= \left| 
\sum_{\sigma \in S_N} \sgn(\sigma) \prod_{n=1}^N 
\left(\alpha_n + \frac{1}{\alpha_n}\right)^{\sigma(n) - 1}
\right|^2,
$$
which we rewrite as
$$
 \sum_{\sigma \in S_N} \sum_{\tau \in S_N} \sgn(\sigma) \sgn(\tau)
\prod_{n=1}^N \left(\alpha_n + \frac{1}{\alpha_n}\right)^{\sigma(n) -
  1}
\left(\overline{\alpha}_n +
\frac{1}{\overline{\alpha}_n}\right)^{\tau(n) - 1}.
$$
Substituting this expression into equation \ref{H_N(s) after
jacobian}, we can write $H_N(s)$ as
\begin{eqnarray*}
&& \frac{1}{2^N N!}\int_{(\C^{\times})^N}\left( \prod_{n=1}^N
\max\left\{|\alpha_n|, |\alpha_n^{-1}|\right\}^{-2s}
\left|\left(\frac{\alpha_n^2 - 1}{\alpha_n^2}
\right)\right|^2 \right) \nonumber  \\ 
&& \hspace{.5cm} \times 
\left( \sum_{\sigma \in S_N} \sum_{\tau \in S_N} \sgn(\sigma) \sgn(\tau)
\prod_{n=1}^N \left(\alpha_n + \frac{1}{\alpha_n}\right)^{\sigma(n) -
  1}
\left(\overline{\alpha}_n +
\frac{1}{\overline{\alpha}_n}\right)^{\tau(n) - 1} \right)
d\lambda_{2N}(\boldsymbol{\alpha}).
\end{eqnarray*}
Exchanging the sums and the integral, and consolidating the products, we
find
\begin{eqnarray*}
&&H_N(s) = \sum_{\sigma \in S_N} \sum_{\tau \in S_N} \sgn(\sigma)
  \sgn(\tau) \frac{1}{2^N N!}\int_{(\C^{\times})^N} 
\left\{\prod_{n=1}^N \max \left\{ |\alpha_n|, \left|
  \alpha_n^{-1} \right| \right\}^{-2s} \right .\\
&& \hspace{1.5cm} \times \left. \left(\frac{\alpha_n^2 -
    1}{\alpha_n^2}\right)\left(\frac{\overline{\alpha}_n^2 -
    1}{\overline{\alpha}_n^2}\right)
\left(\alpha_n + \frac{1}{\alpha_n}\right)^{\sigma(n) - 1}
\left(\overline{\alpha}_n +
\frac{1}{\overline{\alpha}_n}\right)^{\tau(n) - 1} \right\}
d\lambda_{2N}(\boldsymbol{\alpha}).
\end{eqnarray*}
By an application of Fubini's Theorem we find
\begin{equation}
\label{H_N(s) as determinant}
H_N(s) = \frac{1}{N!} \sum_{\sigma \in S_N} \sum_{\tau \in S_N}
\sgn(\sigma) \sgn(\tau) \prod_{n=1}^N \mathcal{I}(\sigma(n), \tau(n))
\end{equation}
where $\mathcal{I}(J,K)$ is given by
$$
\frac{1}{2}
\int_{\C^{\times}} \max \left\{ |\alpha|, \left| \alpha^{-1} \right|
\right\}^{-2s}  
\left(\alpha - \frac{1}{\alpha}\right)
\left(\overline{\alpha} - \frac{1}{\overline{\alpha}}\right) 
\left(\alpha + \frac{1}{\alpha}\right)^{J - 1}
\left(\overline{\alpha} + \frac{1}{\overline{\alpha}}\right)^{K-1} 
\frac{d\lambda_2(\alpha)}{|\alpha|^2}.
$$
Applying Lemma \ref{determinant lemma} to equation \ref{H_N(s) as
determinant} we find $H_N(s)$ is the determinant of the $N \times N$
matrix $\mathcal{I} = \mathcal{I}(J,K)$. 
\end{section}

\begin{section}{The entries of $\mathcal{I}$ are rational functions of
$s$} 

We shall view $\mathcal{I}(J,K)$, not only as an entry in a matrix,
but also a function of $s$.  We note that
$\lambda_2(\alpha)/|\alpha|^2$ is normalized Haar measure on
$\C^{\times}$.  Thus $\mathcal{I}(J,K;s)$ is invariant under the
substitution $\alpha \mapsto \alpha^{-1}$, and we may write
$$
\mathcal{I}(J,K;s) = \int_{\C \setminus D} |\alpha|^{-2s} 
\left(\alpha - \frac{1}{\alpha}\right)
\left(\overline{\alpha} - \frac{1}{\overline{\alpha}}\right)
\left(\alpha + \frac{1}{\alpha}\right)^{J - 1}
\left(\overline{\alpha} + \frac{1}{\overline{\alpha}}\right)^{K-1} 
\frac{d\lambda_2(\alpha)}{|\alpha|^2},
$$ 
where $D$ is the open unit disk.  By setting $\alpha = r e^{i
\theta}$ we may write $\mathcal{I}(J,K;s) = \widehat{h}(J,K;r)$,
where $h(J,K;r)$ is given by 
\begin{equation}
\label{h(J,K;r)}
\int_{0}^{2 \pi} \left(
r e^{i \theta} - \frac{1}{r e^{i \theta}} \right) \left(\frac{r}{e^{i
 \theta}} - \frac{e^{i \theta}}{r} \right)\left(r e^{i \theta} +
\frac{1}{r e^{i \theta}}\right)^{J-1} \left(\frac{r}{e^{i \theta}} +
\frac{e^{i \theta}}{r} \right)^{K-1} \, d\theta
\end{equation}
for $r \in [1, \infty)$, and identically zero on $[0,1)$.

By the change of variables $\theta \mapsto -\theta$ we see that
$h(J,K;r) = h(K,J;r)$. We conclude that $\mathcal{I}$ is a symmetric
matrix whose $J,K$ entry is $\widehat{h}(J,K;s)$.  

\begin{lemma}
\label{rational function lemma}
$\mathcal{I}(J,K;s)$ analytically continues to a rational function.
Specifically 
\begin{equation*}
\mathcal{I}(J,K;s) = \pi \sum_{n = 1}^{N} c_{n}(J) c_{n}(K)
\frac{2s}{s^2 - {n}^2},
\end{equation*}
where 
\begin{equation*}
\label{residues}
c_{n}(J) = \left\{
\begin{array}
{ll}
\left[{J-1 \choose \frac{J+n}{2}}-{J-1 \choose \frac{J+n}{2}
    -1}\right] & \quad \mbox{if} \quad n < J \quad \mbox{and} \quad n
\equiv J \pmod{2} \\
1 & \quad \mbox{if} \quad n=J \\
0 & \quad \mbox{otherwise}
\end{array}
\right. .
\end{equation*}
\end{lemma}
\begin{proof}
Without loss of generality we assume $K \geq J$.  There is a constant
$\mathcal{C}$ such that $h(J,K;r) < \mathcal{C} r^{J+K}$ on $[1,
\infty)$. Thus the integral defining $\mathcal{I}(J,K;s)$ converges
in the half plane $\Re(s) > (J + K)/2$.

Writing $(re^{i \theta} + \frac{1}{re^{i \theta}})^{J-1}$ and
$(\frac{r}{e^{i \theta}} + \frac{e^{i \theta}}{r})^{K-1}$ as sums with
binomial coefficients we may rewrite \ref{h(J,K;r)} as
\begin{eqnarray}
\label{binomial sum}
h(J,K;r) &=& \sum_{j=0}^{J-1} \sum_{k=0}^{K-1} {J-1 \choose j}{K-1
  \choose k} r^{J + K - 2(j + k) -2}  \\
 &\times& \int_0^{2 \pi} \left( r^2 +
\frac{1}{r^{2}} - \left(e^{2 i \theta} + e^{-2 i \theta}\right)
\right) e^{(J-K -2(j-k))i \theta}\, d\theta . \nonumber
\end{eqnarray}
The integral appearing in this expression can be readily evaluated:
\begin{eqnarray}
&& \int_0^{2 \pi} \left(r^2 + \frac{1}{r^2} - \left(e^{2 i
    \theta} + e^{-2 i \theta} \right) \right)e^{(J - K - 2(j-k) ) i
  \theta} \, d\theta \nonumber \\
&& \label{integral cases} \hspace{2cm}= \left\{ 
\begin{array}
{ll}
2\pi \left( r^2 + \frac{1}{r^2} \right) & \quad  k = j+ (K - J)/2 \\
-2\pi & \quad k = j + 1 + (K - J)/2 \\
-2\pi & \quad k = j - 1 + (K - J)/2 
\end{array}
\right. .
\end{eqnarray}
If $J \not\equiv K \pmod 2$ we see that $h(J,K;r)$ (and hence
$\mathcal{I}(J,K;s)$) is identically zero.

The conditions given in equation \ref{integral cases} allow us to
eliminate one of the summations in equation \ref{binomial sum}.  We
use the facts that $0 \leq k \leq K-1$ and $0 \leq j \leq J-1$
together with the conditions in \ref{integral cases} to find
conditions on $j$.  Specifically,
\begin{eqnarray*}
k = j + \frac{K-J}{2} \quad &\Rightarrow& \quad 0 \leq j \leq J-1, \\
k = j + 1 + \frac{K-J}{2} \quad &\Rightarrow& \quad 
0 \leq j \leq
\min\left\{\frac{J+K}{2} - 2, J-1 \right\},
\\
k = j - 1 + \frac{K-J}{2} \quad &\Rightarrow& \quad 
\max\left\{\frac{J-K}{2} + 1, 0 \right\} \leq j \leq J-1.
\end{eqnarray*}
Since $K \geq J$, we can write
$$
\max\left\{\frac{J-K}{2} + 1, 0 \right\} = \delta_{JK}
\quad
\mbox{and}
\quad
\min\left\{\frac{J+K}{2} - 2, J-1 \right\} = J - 1 -\delta_{JK},
$$
where $\delta_{JK} = 1$ if $J=K$ and is 0 otherwise.  From this
information we may write $h(J,K;r)$ as
\begin{eqnarray*}
&& 2\pi \left( \sum_{j=0}^{J-1}{J-1 \choose j}{K-1 \choose
  \frac{K-J}{2} + j}r^{2J-4j} + \sum_{j=0}^{J-1}{J-1 \choose
  j}{K-1 \choose \frac{K-J}{2} + j}r^{2J-4j-4} \right. \\ 
&& \hspace{-.5cm} \left.  -  \sum_{j=\delta_{JK}}^{J-1}{J-1 \choose
  j}{K-1 \choose   \frac{K-J}{2} + j - 1}r^{2J-4j}  
- \sum_{j=0}^{J-1-\delta_{JK}}{J-1
  \choose j}{K-1 \choose \frac{K-J}{2} + j + 1}r^{2J-4j-4}\right).
\end{eqnarray*}
Using the convention that ${K-1 \choose K} = 0$ and ${K-1 \choose -1}
= 0$ we may eliminate $\delta_{JK}$ from the latter two
sums. Reindexing each sum based on the powers of $r$ which appear and
simplifying the binomial coefficients, we find
\begin{eqnarray*}
h(J,K;r) &=& 2\pi \left( \sum_{l = -\frac{J}{2} + 1}^{\frac{J}{2}} {J-1
  \choose \frac{J}{2} - l}{K-1 \choose \frac{K}{2} -l}r^{4l} + 
  \sum_{l = -\frac{J}{2}}^{\frac{J}{2}-1} {J-1 \choose \frac{J}{2} +
  l}{K-1 \choose \frac{K}{2} +l}r^{4l} \right. \\ &-& \left. \sum_{l =
  -\frac{J}{2} + 1}^{\frac{J}{2}} {J-1 \choose \frac{J}{2}
  - l}{K-1 \choose \frac{K}{2} -l -1}r^{4l} - \sum_{l =
  -\frac{J}{2}}^{\frac{J}{2} - 1} {J-1 \choose
  \frac{J}{2} + l}{K-1 \choose \frac{K}{2} +l -1}r^{4l} \right).
\end{eqnarray*}
Note that in the case that $J$ is odd, these sums run over consecutive
odd multiples of $1/2$.  Reindexing the first and third sum by $l
\mapsto -l$, we may combine the first and second sums, and the third
and fourth sums.  We may then write $h(J,K;r)$ as
$$
 2\pi \left(\sum_{l = -\frac{J}{2}}^{\frac{J}{2}-1} {J-1 \choose \frac{J}{2} +
  l}{K-1 \choose \frac{K}{2} +l}\left(r^{4l} + r^{-4l}\right)
- \sum_{l = -\frac{J}{2}}^{\frac{J}{2} - 1} {J-1
   \choose \frac{J}{2} + l}{K-1 \choose \frac{K}{2} +l
   -1}\left(r^{4l} + r^{-4l}\right) \right).
$$ 
Due to the symmetry in the summands we may reindex the sums using
only positive indices.  Let $l_0 = 0$ if $J$ and $K$ are even,
and $l_0 = 1/2$ if $J$ and $K$ are odd, then
\begin{eqnarray*}
h(J,K;r) &=&
2\pi\left[
{K-1 \choose \frac{J+K}{2} - 1} - {K-1 \choose \frac{J+K}{2}}
\right]\left(r^{2J} + r^{-2J}\right) \\
&+& 2\pi \sum_{l = l_0}^{\frac{J}{2} -1} \left[{J-1 \choose
    \frac{J}{2} + l} - {J-1 \choose \frac{J}{2} + l -1}
\right]
\left[
{K-1 \choose \frac{K}{2} + l} - {K-1 \choose \frac{K}{2} + l -1}
\right]\left(r^{4l} + r^{-4l} \right).
\end{eqnarray*}

We are now in position to compute $\widehat{h}(J,K;s)$.  There is a
correspondence between the coefficients and powers of $r$ which appear
in $h(J,K;r)$ and the poles and residues of $\widehat{h}(J,K;s)$.  As
was demonstrated in the proof of Corollary \ref{mellin inversion
corollary}, the Mellin transform of $r^{4l} + r^{-4l}$ analytically
continues to the rational function $s/(s^2 - 4 l^2)$.  Thus
$\mathcal{I}(J,K;s)$ extends to a rational function: 
\begin{eqnarray*}
\mathcal{I}(J,K;s) &=& 2\pi\left[
{K-1 \choose \frac{J+K}{2}-1} - {K-1 \choose \frac{J+K}{2}}
\right] \frac{s}{s^2 - J^2} \\
&+& 2\pi \sum_{l = l_0}^{\frac{J}{2} -1} \left[{J-1 \choose
    \frac{J}{2} + l} - {J-1 \choose \frac{J}{2} + l -1}
\right]
\left[
{K-1 \choose \frac{K}{2} + l} - {K-1 \choose \frac{K}{2} + l -1}
\right]\frac{s}{s^2 - 4 l^2}.
\end{eqnarray*}
Reindexing this sum by setting $n = 2l$ we find 
\begin{eqnarray*}
\mathcal{I}(J,K;s) &=& 2\pi\left[
{K-1 \choose \frac{J+K}{2}-1} - {K-1 \choose \frac{J+K}{2}}
\right] \frac{s}{s^2 - J^2} \\
&+& 2\pi \sum_n \left[{J-1 \choose
    \frac{J+n}{2}} - {J-1 \choose \frac{J+n}{2}-1}
\right]
\left[
{K-1 \choose \frac{K+n}{2}} - {K-1 \choose \frac{K+n}{2}-1}
\right]\frac{s}{s^2 - n^2},
\end{eqnarray*}
where the sum is over $n \in \{1,3,\ldots, J-2\}$ if $J$ and $K$ are
odd, and over $n \in \{2,4,\ldots,J-2\}$ if $J$ and $K$ are even.
If $J = K$ the leading coefficient is $1$.  If we set
\begin{equation*}
\label{uniform coefficient}
c_{n}(J) = \left\{
\begin{array}
{ll}
\left[{J-1 \choose \frac{J+n}{2}}-{J-1 \choose \frac{J+n}{2}
    -1}\right] & \quad \mbox{if} \quad n < J \quad \mbox{and} \quad n
\equiv J \pmod{2} \\
1 & \quad \mbox{if} \quad n=J \\
0 & \quad \mbox{otherwise}
\end{array}
\right. ,
\end{equation*}
then we may write 
\begin{equation}
\label{uniform formulation}
\mathcal{I}(J,K;s) = 
\pi \sum_{n = 1}^{N} c_{n}(J) c_{n}(K)
\frac{2s}{s^2 - {n}^2}.
\end{equation}
It is easy to verify that this expression is symmetric in $J$ and $K$,
giving $\mathcal{I}(J,K;s) = \mathcal{I}(K,J;s)$ as expected.
Additionally if $J \not \equiv K \pmod 2$, the expression in equation
\ref{uniform formulation} yields $\mathcal{I}(J,K;s) = 0$.  This
proves the lemma. 
\end{proof}

We identify $\mathcal{I}(J,K;s)$ with the rational function it extends
to.  When $J$ and $K$ are odd, $\mathcal{I}(J,K;s)$ has poles at $\pm
1, \pm 3,\ldots, \pm \min\{J,K\}$.  When $J$ and $K$ are even,
$\mathcal{I}(J,K;s)$ has poles at $\pm 2$, $\pm 4$, $\ldots$, $\pm
\min\{J,K\}$.  $\mathcal{I}(J,K;s)$ has a zero of multiplicity one at
$0$.   

We are now in position to prove the first part of Theorem \ref{main
theorem}.  $H_N(s)$ is the determinant of $\mathcal{I}$, and the
entries of $\mathcal{I}$ extend to rational functions of $s$. Since the 
determinant is a polynomial in the entries of a matrix,
$H_N(s)$ itself extends to a rational function of $s$.  In fact, since
the determinant is a homogeneous polynomial in the entries of a matrix
and the entries of $\mathcal{I}$ analytically continue to odd
functions, $H_N(s)$ analytically continues to an even rational
function when $N$ is even, and analytically continues to an odd 
rational function when $N$ is odd.  We also see $H_N(s)$ has a zero of
multiplicity $N$ at $0$.   
\end{section}

\begin{section}{$H_N(s)$ is a simple product}
In this section we express $\det(\mathcal{I})$ as a simple product.
The structure of the poles and residues of $\mathcal{I}(J,K;s)$
will allow us to find linear dependence relations on the rows of
$\mathcal{I}$.  

Let $B_{n}$ be the $N \times N$ matrix whose $J,K$ entry is
the integer $c_{n}(J)c_{n}(K)$. Then by Lemma
\ref{rational function lemma} we have the matrix equation
$$
\mathcal{I} = \sum_{n = 1}^{N} B_{n} \, \frac{2 \pi s}{s^2 -
  n^2}.
$$
Define $\boldsymbol{\omega}^{\T}_{n} \in \Q^N$ to be the
row vector given by $\boldsymbol{\omega}^{\T}_{n} = ( c_{n}(K)
)_{K=1}^N$.  It follows then that the $J$th row vector of
$B_{n}$ is given by $c_{n}(J) \boldsymbol{\omega}^{\T}_{n}$, and thus
every row of $B_{n}$ is a scalar multiple of
$\boldsymbol{\omega}^{\T}_{n}$.  

We may find a nonzero vector $\boldsymbol{\psi} \in \Q^N$ such that $\ 
\boldsymbol{\omega}_{n}^{\T}\boldsymbol{\psi} = 0$ for $1 \leq
n \leq N-1$.  In fact, $B_{n}\boldsymbol{\psi} =
\boldsymbol{0}$ for $1 \leq n \leq N-1$, leading us to the
vector equation 
$$
\mathcal{I}\boldsymbol{\psi} = \sum_{n=1}^N
B_{n}\boldsymbol{\psi} \frac{2 \pi s}{s^2 - n^2} = 
B_N\boldsymbol{\psi} \frac{2\pi s}{s^2 - N^2}.
$$
We see that $ \left(\mathcal{I} - B_N \frac{2\pi s}{s^2 -
  N^2}\right)\boldsymbol{\psi} = \boldsymbol{0},$ and so $ \det
\left(\mathcal{I} - B_N \frac{2\pi s}{s^2 -  N^2}\right) = 0.$
From the definition of $B_N$ and Lemma \ref{rational function lemma}, we find
$$
\mathcal{I} - B_N \frac{2\pi s}{s^2 - N^2} = 
\begin{pmatrix}
\mathcal{I}(1,1)  & \mathcal{I}(1,2) & \ldots &  \mathcal{I}(1,N)\\
\mathcal{I}(2,1)  & \mathcal{I}(2,2) & \ldots &  \mathcal{I}(2,N)\\
\vdots    & \vdots    & \ddots & \vdots    \\
\mathcal{I}(N,1)  & \mathcal{I}(N,2) & \ldots &  \mathcal{I}(N,N) -
\frac{2 \pi s}{s^2 - N^2} \\
\end{pmatrix}.
$$
Taking determinants and exploiting the multilinearity of the determinant
$$
\det \begin{pmatrix}
\mathcal{I}(1,1)  & \mathcal{I}(1,2) & \ldots &  \mathcal{I}(1,N)\\
\mathcal{I}(2,1)  & \mathcal{I}(2,2) & \ldots &  \mathcal{I}(2,N)\\
\vdots    & \vdots    & \ddots & \vdots    \\
\mathcal{I}(N,1)  & \mathcal{I}(N,2) & \ldots &  \mathcal{I}(N,N)\\
\end{pmatrix}
=
\det \begin{pmatrix}
\mathcal{I}(1,1)  & \mathcal{I}(1,2) & \ldots &  0\\
\mathcal{I}(2,1)  & \mathcal{I}(2,2) & \ldots &  0\\
\vdots    & \vdots    & \ddots & \vdots    \\
\mathcal{I}(N,1)  & \mathcal{I}(N,2) & \ldots &  \frac{2 \pi s}{s^2 - N^2} \\
\end{pmatrix}.
$$
The left hand side is $H_N(s)$.  By a simple induction argument we
finally arrive at a simple product formulation of $H_N(s)$:
$$H_N(s) = \prod_{n=1}^N \frac{2 \pi s}{s^2 - n^2}.$$

\end{section}

\bibliography{recstarbody3}
\end{document}